\documentclass[a4paper]{amsart}

\usepackage[utf8]{inputenc}
\usepackage[T1]{fontenc}
\usepackage{lmodern}
\usepackage{amssymb}

\usepackage[pdfusetitle]{hyperref}

\let\newterm\emph

\def\ie{i.\,e.}
\def\cf{cf.}

\def\arxiv#1{\href{http://arxiv.org/abs/#1}{\texttt{arXiv:#1}}}
\def\doi#1{\href{http://dx.doi.org/#1}{doi:#1}}

\def\Z{{\mathbb Z}}
\def\Q{{\mathbb Q}}
\def\C{{\mathbb C}}

\def\m{{\mathfrak m}}
\def\p{{\mathfrak p}}
\def\t{{\mathfrak t}}

\DeclareMathOperator{\Ann}{Ann}

\DeclareMathOperator{\depth}{depth}
\DeclareMathOperator{\Ext}{Ext}
\DeclareMathOperator{\Supp}{Supp}

\theoremstyle{plain}
\newtheorem{theorem}{Theorem}[section]
\newtheorem{proposition}[theorem]{Proposition}
\newtheorem{lemma}[theorem]{Lemma}
\newtheorem{corollary}[theorem]{Corollary}
\newtheorem{claim}[theorem]{Claim}

\theoremstyle{definition}

\theoremstyle{remark}
\newtheorem{remark}[theorem]{Remark}
\newtheorem*{acknowledgements}{Acknowledgements}

\numberwithin{equation}{section}

\def\CM{Cohen--Macaulay}
\def\CS{Chang--Skjelbred}

\begin{document}

\title{Exact sequences for equivariantly formal spaces}
\author{Matthias Franz and Volker Puppe}
\thanks{The first author was partially supported by
  grants from DAAD and NSERC}
\address{Department of Mathematics, University of Western Ontario,
  London, Ont.\ N6A\;5B7, Canada}
\email{mfranz@uwo.ca}
\address{FB Mathematik und Statistik, Universität Konstanz,
  78457 Konstanz, Germany}
\email{volker.puppe@uni-konstanz.de}

\subjclass[2010]{Primary 55N91; secondary 13C14}

\begin{abstract}
  Let $T$ be a torus.
  We present an exact sequence
  relating the relative equivariant cohomologies of the skeletons
  of an equivariantly formal $T$-space.
  This sequence, which goes back to Atiyah and Bredon,
  generalizes the so-called \CS~lemma.
  As coefficients, we allow prime fields and subrings of the rationals,
  including the integers.
  We extend to the same coefficients
  a generalization of this ``Atiyah--Bredon sequence''
  for actions without fixed points
  which has recently been obtained by Goertsches and Töben.
\end{abstract}

\maketitle

\section{Introduction}

Let
$T=(S^1)^n$ be a torus of dimension~$n$ and $X$ a ``nice''
$T$-space.
Suppose that $X$ is
\newterm{cohomologically equivariantly formal}%
\footnote{It seems that the terminology `equivariantly formal'
was coined by~\cite{GoreskyKottwitzMacPherson:98}. We add
`cohomologically' to avoid confusion with notions in
rational homotopy theory.} (or \newterm{CEF}),
which means that its equivariant cohomology~$H_T^*(X)$
(with rational coefficients)
is free over the polynomial ring~$H^*(BT)$.
Then the so-called {\CS} lemma~\cite[(2.3)]{ChangSkjelbred:74}
asserts that the sequence
\begin{equation}
  \label{C-S}
  \let\longrightarrow\rightarrow
  0\longrightarrow
  H_T^*(X)\longrightarrow
  H_T^*(X^T)\longrightarrow
  H_T^{*+1}(X_1,X^T)
\end{equation}
is exact, where $X^T\subset X$ denotes the fixed point set
and $X_1$ the union of all orbits of dimension at most~$1$.
In other words, $H_T^*(X)$ coincides, as subalgebra of~$H_T^*(X^T)$,
with the image of~$H_T^*(X_1)\to H_T^*(X^T)$.
Usually $X^T$~and~$X_1$ are much simpler than~$X$,
so that \eqref{C-S} gives an easy method
to compute $H_T^*(X)$. This can be used for instance for a short proof of
Jurkiewicz's theorem about the cohomology of smooth projective
toric varieties because they are known
to be equivariantly formal, see Remark~\ref{conditions-cef}.
There are many other applications of the {\CS} lemma in algebraic and
symplectic geometry, often with implications for combinatorics.

Roughly at the same time as Chang and Skjelbred, Atiyah proved a much more
general theorem in the context of
equivariant $K$-theory~\cite[Ch.~7]{Atiyah:74}.
Bredon~\cite{Bredon:74} then observed that it applies equally to cohomology.
In the context of toric varieties an Atiyah-like theorem
was proven
by Barthel--Brasselet--Fieseler--Kaup
for cohomology and intersection homology
\cite[Thm.~4.3]{BarthelBrasseletFieselerKaup:02}.
Recently
Goertsches and Töben~\cite{GoertschesToeben:??} made
the very interesting observation that the {\CM} property,
which Atiyah used in an essential way,
gives generalizations also for fixed-point-free actions.

The principal aim of this note is to present a generalized
cohomological version of Atiyah's theorem.
Our proof in Section~\ref{atiyah} will be a modification of his,
adapted to cohomology and with a special emphasis
on coefficient rings~$R$ other than the rationals,
namely on subrings of~$\Q$ and on prime fields.
As a corollary, we get an integral version of the Chang--Skjelbred lemma
with certain restrictions on the isotropy groups.
We also extend the Goertsches--Töben result to the same coefficient
rings and to more general spaces than compact $T$-manifolds.

\begin{acknowledgements}
  The authors thank Winfried Bruns, Jean-Claude Hausmann,
  Lex Renner and Sue Tolman for helpful discussions.
  We are also indebted to the anonymous referee
  for suggesting improvements to the presentation of the paper.
\end{acknowledgements}

\section{The Atiyah--Bredon sequence}\label{atiyah}

We use Alexander--Spanier cohomology,
\cf~\cite[Ch.~8]{Massey:78}.
Coefficients are taken in a subring~$R$ of~$\Q$
or in the finite field~$R=\Z_p$.
Here the letter~$p$ denotes a prime, as it does in the rest of the paper.
All topological spaces will be assumed to be paracompact and Hausdorff.
Recall that on locally contractible spaces
Alexander--Spanier cohomology and singular cohomology coincide.

Let $T=(S^1)^n$ be a torus and $X$ a $T$-space.
(Note that its Borel construction~$X_T$ is again paracompact Hausdorff.)
We assume $X$ to be finite-dimensional,
with only finitely many orbit types
and such that $H^*(X^{K})$ is finitely generated over~$R$
for all closed subgroups~$K$ of~$T$.%
\footnote{These assumptions can be weakened,
\cf~\cite[Sec.~3.2]{AlldayPuppe:93}. For example,
$X$ can be finitistic instead of finite-dimensional, and, for~$R=\Q$,
must only have finitely many \emph{connected} orbit types.}
For example, $X$ may be a finite $T$-CW~complex,
a compact topological $T$-manifold
or a complex algebraic variety
with an algebraic action of
~$(\C^\times)^n\supset T$.

Denote by~$X_i$, $-1\le i\le n$, the equivariant $i$-skeleton of~$X$, \ie,
the union of all orbits of dimension~$\le i$.
In particular, $X_{-1}=\emptyset$, $X_0=X^T$ and $X_n=X$.
Since $X$ is Hausdorff with only finitely many orbit types,
each~$X_i$ is closed in~$X$.
The maximal $p$-torus~$T_p\subset T$
is the subgroup of all elements of order~$p$
(together with~$1\in T$). It is isomorphic to~$\Z_p^n$.
The $T_p$-equivariant $i$-skeleton~$X_{p,i}$ of~$X$ is
the set of all $T_{p}$-orbits
consisting of at most $p^i$~points.
One clearly has $X_i\subset X_{p,i}$.

The inclusion of pairs~$(X_i,X_{i-1})\hookrightarrow(X,X_{i-1})$
gives rise to a long exact sequence
\begin{equation}\label{exact-triple}
  \let\longrightarrow\rightarrow
  \cdots
  \longrightarrow H^*_T(X, X_i)
  \longrightarrow H^*_T(X, X_{i-1})
  \longrightarrow 
  \let\longrightarrow\rightarrow
  H^*_T(X_i, X_{i-1})
  \stackrel\delta\longrightarrow H^{*+1}_T(X, X_i)
  \longrightarrow \cdots,
\end{equation}
and likewise $(X_{i+1},X_{i-1})\hookrightarrow(X_{i+1},X_i)$ gives maps
\begin{equation}
  H^*_T(X_i, X_{i-1})\to H^{*+1}_T(X_{i+1}, X_i).
\end{equation}

\begin{theorem}
  \label{Atiyah-Bredon}
  Let $X$ be a $T$-space such that $H^*_T(X)$ is a free $H^*(BT)$-module.
  Suppose in addition that for all~$i$
  \begin{subequations}\label{condition-skeleton}
  \begin{align}
        X_{p,i} &= X_i \label{condition-skeleton-field}\\
    \intertext{if $R=\Z_p$, or, in case~$R\subset\Q$,}
        X_{p,i-1} &\subset X_i \label{condition-skeleton-subring}
  \end{align}
  \end{subequations}
  for all~$p$ not invertible in~$R$.
  Then the following ``Atiyah--Bredon sequence'' is exact:
  \begin{equation}\label{sequence-atiyah}
  \let\longrightarrow\rightarrow
    0
    \longrightarrow H^*_T(X)
    \longrightarrow H^*_T(X_0)
    \longrightarrow H^{*+1}_T(X_1, X_0)
    \longrightarrow \cdots
    \longrightarrow H^{*+n}_T(X_n, X_{n-1})
    \longrightarrow 0.
  \end{equation}
  More generally, if condition
  \eqref{condition-skeleton-field}~or~\eqref{condition-skeleton-subring}
  holds for all~$i\le k$,
  then one still has exactness for 
  \begin{equation}\label{sequence-atiyah-k}
  \let\longrightarrow\rightarrow
    0
    \longrightarrow H^*_T(X)
    \longrightarrow H^*_T(X_0)
    \longrightarrow H^{*+1}_T(X_1, X_0)
    \longrightarrow \cdots
    \longrightarrow H^{*+k}_T(X_k, X_{k-1}).
  \end{equation}
\end{theorem}

The content of the {\CS} lemma is just exactness
of the Atiyah--Bredon sequence at the first terms.
Rephrasing Theorem~\ref{Atiyah-Bredon} for~$k\in\{0,1\}$ and
integer coefficients, we get:

\begin{corollary}\label{Chang-Skjelbred-integers}
  Assume $R=\Z$ and let $X$ be a $T$-space such that $H_T^*(X)$ is free over~$H^*(BT)$.
  Then
  $$
  \let\longrightarrow\rightarrow
    0\longrightarrow
    H_T^*(X)\longrightarrow
    H_T^*(X_0)
  $$
  is exact. If in addition the isotropy group of each~$x\not\in X_1$
  is contained in a proper subtorus of~$T$
  (\ie, if condition (2.3b) holds for~$i\le 1$),
  then the following sequence is exact:
  $$
  \let\longrightarrow\rightarrow
    0\longrightarrow
    H_T^*(X)\longrightarrow
    H_T^*(X_0)\longrightarrow
    H_T^{*+1}(X_1,X_0).
  $$
\end{corollary}

The injectivity of the map~$H_T^*(X)\to H_T^*(X_0)$
is an immediate consequence of the localization theorem
in equivariant cohomology, see Lemma~\ref{first-claim} below.

The above ``proper subtorus condition''
is violated in the counterexample given
by Tolman--Weitsman~\cite[Sec.~4]{TolmanWeitsman:99}
because there all orbits of dimension~$2$ have
isotropy group~$\Z_2\times\Z_2\subset S^1\times S^1$.
The full condition~\eqref{condition-skeleton-subring}
can be stated analogously: For~$R=\Z$
the isotropy group~$T_x$ of any~$x\in X$ must be contained
in a subtorus of dimension $\le\dim T_x+1$.

\begin{remark}\label{conditions-cef}
  In some cases general criteria guarantee
  the freeness of~$H_T^*(X)$ over $H^*(BT)$, for instance
  if $H^*(X)$ is free over~$R$ and concentrated in even degrees
  (Leray--Hirsch).
  Assume that $H^*(X)$ is free over~$R\subset\Q$. Then
  $H_T^*(X)$ is free over~$H^*(BT)$
  if $X$ is a compact Kähler manifold with~$X^T\ne\emptyset$
  (for example, a smooth projective $(\C^\times)^n$-variety)
  \cite[Thm.~II.1.2]{Blanchard:56} 
  or
  if $X$ is a compact Hamiltonian $T$-manifold
  \cite{Frankel:59}. 
  If all isotropy groups of~$X$ are connected, then
  the exactness of the Atiyah--Bredon sequence with integer
  coefficients can be characterized by a homological condition
  on $H_T^*(X)$, see~\cite{FranzPuppe:07}.
\end{remark}

\section{Isotropy groups}

Before proving Theorem~\ref{Atiyah-Bredon},
we want to reformulate the conditions
\eqref{condition-skeleton-field}~and \eqref{condition-skeleton-subring}
on the equivariant skeletons in algebraic terms.
In the sequel, we will refer to
\eqref{condition-skeleton-field}~and~\eqref{condition-skeleton-subring}
together as condition~\eqref{condition-skeleton}.
Notice that this condition is always satisfied for~$R=\Q$.

For any closed subgroup~$T'\subset T$
there is an isomorphism~$T\cong(S^1)^n$
and a divisor chain~$m_q \mathbin| m_{q-1}\mathbin|\cdots\mathbin|m_1$
such that $T'$ corresponds to the subgroup
\begin{equation}
  \label{decomp-subtorus}
  \Z_{m_1}\times\cdots\times\Z_{m_q}\times(S^1)^r\subset(S^{1})^{n}.
\end{equation}
(To see this, consider the inverse image of~$T'$ under
the exponential map~$\exp\colon\t\to T$
and compare it with the lattice~$\exp^{-1}(1)$.)

Recall that the dimension of a 
module~$M$
over a ring~$S$ is the Krull dimension of~$S/\Ann M$,
that is, the length of a maximal chain
$\p_0\subsetneq\p_1\subsetneq\cdots$ 
of prime ideals in~$S$,
all containing the annihilator~$\Ann M$ of~$M$.
We denote the dimension of~$R$ over itself by~$d$.
Hence, $d=0$ if $R$ is a field, and $d=1$ if~$R\subsetneq\Q$.
Since $R$ is a PID, the dimension
of a polynomial ring~$R[t_{1},\dots,t_{s}]$ is $d+s$
\cite[Prop.~III.13]{Serre:65}.

\begin{lemma}\label{cohomological-dimension-isotropy}
  Let $T'\subset T$ be a closed subgroup and
  let $s$ be the number of~$m_j$'s
  in the decomposition~\eqref{decomp-subtorus}
  which are not invertible in~$R$.
  Then
  $$
    \dim_{H^*(BT)}H^*(BT')=
    \begin{cases}
      r+s & \hbox{if $R$ is a field or $s>0$,} \\
      r+1 & \hbox{if $R\subsetneq\Q$ and $s=0$.}
    \end{cases}
  $$
\end{lemma}

\begin{proof}
  The $H^*(BT)$-action on~$H^*(BT')$ comes from the algebra map
  induced by the inclusion~$T'\hookrightarrow T$.
  Hence, the dimension of the $H^*(BT)$-module~$H^*(BT')$
  is equal to the Krull dimension of the (evenly graded) image~$A$
  of~$H^*(BT)$ in~$H^*(BT')$.

  We choose isomorphisms
  $H^{2*}(BS^1)\cong R[t]$ and $H^{2*}(B\Z_m)\cong R[t]/m t$
  such that the map~$H^{2*}(BS^1)\to H^{2*}(B\Z_m)$ corresponds
  to the canonical quotient map~$R[t]\to R[t]/m t$.
  By the Künneth theorem for cohomology, we find
  \begin{align*}
    A
    &\cong H^{2*}(B\Z_{m_1})\otimes\cdots\otimes H^{2*}(B\Z_q)\otimes
      H^{2*}(BS^1)^{\otimes r} \\
    &\cong R[t_1,\ldots,t_q,t'_1,\ldots,t'_r]/(m_1 t_1,\ldots,m_q t_q) \\
    &= R[t_1,\ldots,t_s,t'_1,\ldots,t'_r]/(m_1 t_1,\ldots,m_s t_s)
  \end{align*}
  since $m_{s+1}$,~\ldots,~$m_{q}$ are invertible in~$R$.

  $R$, a field or~$s=0$:
  then $A\cong R[t_1,\ldots,t_s,t'_1,\ldots,t'_r]$,
  so the assertion is clear.

  $R\subsetneq\Q$~and~$s>0$:
  An ascending chain of prime ideals in~$A$ of length~$k$ is the same as an
  ascending chain~$\p_0\subset\cdots\subset\p_k$ of prime ideals
  in~$B=R[t_1,\ldots,t_s,t'_1,\ldots,t'_r]$
  with~$(m_1 t_1,\ldots,m_s t_s)\subset\p_0$.
  Since one can clearly find such a chain of length~$r+s$,
  we have to show that there are no longer chains.
  If $\p_0$ contains some prime factor~$p$ of~$m_s$,
  then $B/\p_0$ is a quotient of~$\Z_p[t_1,\ldots,t_s,t'_1,\ldots,t'_r]$,
  which is of dimension~$r+s$.
  Otherwise, $t_s\in\p_0$
  and $B/\p_0$ is a quotient of~$R[t_1,\ldots,t_{s-1},t'_1,\ldots,t'_r]$,
  which by induction has again dimension~$r+s$.
  In any case, we find that $k\le r+s$.
\end{proof}

\begin{proposition}\label{condition-isotropy}
  The conditions~\eqref{condition-skeleton} for all~$i\le k$
  are equivalent to the conditions
  \begin{equation}\label{condition-isotropy-display}
    \forall x\in X\quad
      \dim_{H^*(BT)}H^*(BT_x)\ge d+n-i \;\Longrightarrow\; x\in X_i
  \end{equation}
  for all~$i\le k$.
\end{proposition}

\begin{proof}
  $R=\Q$:
  Since the dimension of~$H^*(BT_x)$ is just
  the rank of~$T_{x}$, condition~\eqref{condition-isotropy-display}
  holds for all spaces, as does~\eqref{condition-skeleton}.

  $R=\Z_p$: By Lemma~\ref{cohomological-dimension-isotropy},
  condition~\eqref{condition-isotropy-display} means that
  the rank of the maximal $p$-torus contained in~$T_x$ equals
  the dimension of~$T_x$ for~$x\in X_k$, and differs by at most~$i-k-1$
  for~$x\in X_i\setminus X_{i-1}$, $i>k$.
  This is equivalent to~$X_{p,i}=X_i$ for~$i\le k$.

  $R\subsetneq\Q$: Here condition~\eqref{condition-isotropy-display}
  translates into the following: for each non-invertible prime~$p$
  the rank of the maximal $p$-torus in~$T_x$ differs by at most~$1$
  from the dimension of~$T_x$ for~$x\in X_k$, and by at most~$i-k$
  for~$x\in X_i\setminus X_{i-1}$, $i>k$.
  This is the same as saying $X_{p,i-1}\subset X_i$ for all~$i\le k$.
\end{proof}

Note that
conditions~\eqref{condition-isotropy-display} are always
true in $K$-theory (see~\cite[(7.1)]{Atiyah:74}),
essentially because the representation ring
of a finite group is a finitely generated $\Z$-module.

\section{Proof of Theorem~\protect\ref{Atiyah-Bredon}}
\label{proof-Atiyah-Bredon}

Before starting in earnest, we remark that our assumptions on~$X$ imply
that all $R$-modules~$H^{*}(X_{j},X_{i})$, $j\ge i$, are finitely generated,
hence so are the $H^*(BT)$-modules~$H_{T}^{*}(X_{j},X_{i})$.
(Use the usual long exact sequences for the first claim.
For equivariant cohomology, the proof given in~\cite[(3.10.1)]{AlldayPuppe:93}
for~$R=\Q$ still applies.)

\smallskip

  The sequences 
  \eqref{sequence-atiyah}~and~\eqref{sequence-atiyah-k}
  are exact
  if and only if for all~$i\ge0$ (respectively $0\le i\le k$)
  the sequence~\eqref{exact-triple} splits into short exact sequences
  \begin{equation}\label{short-sequence}
  \let\longrightarrow\rightarrow
    0
    \longrightarrow H^*_T(X,X_{i-1})
    \longrightarrow H^*_T(X_i,X_{i-1})
    \longrightarrow H^{*+1}_T(X,X_i)
    \longrightarrow 0,
  \end{equation}
  \ie, if and only if the
  inclusion of pairs~$(X,X_{i-1})\hookrightarrow(X,X_i)$
  induces the zero map in cohomology.
  (See \cite[Lemma~4.1]{FranzPuppe:07}.)

Let $M$ be a graded $H^*(BT)$-module.
Then any prime ideal associated with~$M$
is homogeneous, \ie, generated by its homogeneous elements
\cite[Lemma~1.5.6]{BrunsHerzog:93},
hence contained in some maximal homogeneous ideal~$\m$.
This implies $M=0$ if and only if $M_\m=0$
for all maximal homogeneous ideals~$\m\subset H^*(BT)$.
It suffices therefore to prove exactness after localizing
at such an ideal~$\m$. Note that it always contains $H^{>0}(BT)$.
From now on all localized modules will be over the
regular local ring~$A=H^*(BT)_\m$, whose maximal ideal
we denote by~$\m_{A}$.

By induction on~$i\ge0$ (resp.~$0\le i\le k$), we show:

\begin{claim}\label{big-claim}
  The sequence~\eqref{short-sequence} is exact, and
  $H^*_T(X,X_i)_\m$ is zero or a {\CM} module of dimension~$ d+n-i-1$.
\end{claim}

Recall that the \newterm{depth} of a finitely generated $A$-module~$M$
is the maximal length of an $M$-regular sequence in the maximal
ideal~$\m_{A}$.
An $M$-regular sequence is a sequence~$a_1$,~\ldots,~$a_l\in\m_{A}$
such that $(a_1,\ldots,a_l)M\neq M$ and such that $a_j$ is not
a zero divisor on~$M/(a_1,\ldots,a_{j-1})M$ for all~$j\le l$.
One always has
\begin{equation}\label{depth-le-dim}
  \depth M\le\dim M;
\end{equation}
if $M\ne0$ and
\begin{equation}
  \depth M=\dim M,
\end{equation}
then $M$ is called \newterm{\CM}.
(We use \cite{Serre:65} and \cite{BrunsHerzog:93} as
references for commutative algebra.
The reader might also
find the summary of results 
in~\cite[App.~A]{AlldayPuppe:93} helpful.
They were compiled with applications in equivariant cohomology in mind.)

\medbreak

We start the proof with several lemmas.
The first two are standard.

\begin{lemma}\label{depth-cokernel}
  Let $0\to U\to M\to N\to 0$ be an exact sequence of $A$-modules. Then
  \[
    \depth N\ge\min\bigl\{\depth U-1,\depth M\bigr\}.
  \]
\end{lemma}

\begin{proof}
  Since
  $$
     \depth M=\inf\bigl\{i:\Ext^i(K, M)\ne0\bigr\},
  $$
  where $K$ denotes the residue field of~$A$,
  the assertion follows from the long exact sequence
  \begin{equation*}
\let\longrightarrow\rightarrow
    \cdots
    \longrightarrow\Ext^i(K, U)
    \longrightarrow\Ext^i(K, M)
    \longrightarrow\Ext^i(K, N)
    \longrightarrow\Ext^{i+1}(K, U)
    \longrightarrow\cdots. \qedhere
  \end{equation*}
\end{proof}

\begin{lemma}\label{dimension-submodule}
  Let $M$ be a {\CM} $A$-module.
  If $U\subset M$ is a non-zero submodule, then $\dim U=\dim M$.
\end{lemma}

\begin{proof}
  Choose a prime ideal~$\p$ associated to~$U$, hence also to~$M$.
  Since $M$ is {\CM}, we have
  $\dim A/\p=\dim M$ by~\cite[Prop.~IV.13]{Serre:65}.
  On the other hand, $\dim A/\p\le\dim U\le\dim M$.
\end{proof}

\begin{lemma}\label{first-claim}
  For all~$i\le k$ one has
  $ 
    \dim H^*_T(X,X_i)\le d+n-i-1
  $. 
\end{lemma}


\begin{proof}
  This follows from the localization theorem in equivariant cohomology:
  For a prime ideal~$\p\subset H^*(BT)$, define
  $
    X^\p=\{\,x\in X: \p\supset\Ann H^*(BT_x)\,\}.
  $
  The localization theorem asserts that
  the inclusion~$X^\p\hookrightarrow X$ induces an isomorphism
  \begin{equation*} 
    H^*_T(X)_\p=H^*_T(X^\p)_\p,
  \end{equation*}
  see for example \cite[Sec.~3.2]{AlldayPuppe:93}.
  An analogous formula holds for relative cohomology.
  
  Now assume that $\dim H^*_T(X,X_i)\ge d+n-i$.
  This means that there is a prime ideal~$\p\supset\Ann H^*_T(X,X_i)$
  such that $\dim H^*(BT)/\p\ge d+n-i$.
  Hence $\dim H^*(BT_x)
  \ge d+n-i$ for any~$x\in X^\p$.
  By Proposition~\ref{condition-isotropy}, this implies $x\in X_i$,
  \ie, $X^\p\subset X_i$. The localization theorem now gives
  $$
    H_T^*(X,X_i)_\p=H_T^*(X^\p,X^\p_i)_\p=H_T^*(X^\p,X^\p)_\p=0.
  $$
  But this is impossible because for any prime ideal~$\p$
  one has $\p\supset\Ann H^*_T(X,X_i)$ if and only if $H_T^*(X,X_i)_\p\ne0$
  \cite[Prop.~I.3]{Serre:65}.
  (Here we are using that $H_T^*(X,X_i)$ is finitely generated over~$H^*(BT)$.)
  Hence $\dim H^*_T(X,X_i)\le d+n-i-1$.
\end{proof}

\begin{lemma}\label{second-claim}
  For all~$i\ge0$ one has
  $ 
    \depth H^*_T(X_i,X_{i-1})_\m\ge n-i
  $. 
\end{lemma}

\begin{proof}
  Each~$x\in X_i\setminus X_{i-1}$ is fixed by exactly one
  $(n-i)$-dimensional torus~$T_\alpha\subset T$.
  Because $X$ has only finitely many orbit types,
  we can write $X_i\setminus X_{i-1}$ as the disjoint union of finitely
  many subsets~$Y_\alpha=(X_i\setminus X_{i-1})^{T_{\alpha}}$.
  Since $\bar Y_\alpha\subset Y_\alpha\cup X_{i-1}$,
  a Mayer--Vietoris argument gives
  $$
    H_T^*(X_i,X_{i-1})
    =\bigoplus_\alpha H^*_T(\bar Y_\alpha, \bar Y_\alpha\cap X_{i-1}).
  $$
  (Recall that with respect to Alexander--Spanier cohomology,
  the Mayer--Vietoris sequence is exact for all pairs
  of closed subspaces, \cf~\cite[Sec.~8.6]{Massey:78}.)

  We have
  \begin{equation*} 
    M := H^*_T(\bar Y_\alpha, \bar Y_\alpha\cap X_{i-1})
       = H^*_{T/T_\alpha}(\bar Y_\alpha, \bar Y_\alpha\cap X_{i-1})
           \otimes_R H^*(BT_\alpha)
  \end{equation*}
  because $\bar Y_\alpha$ is fixed by~$T_\alpha$.
  A sequence of generators~$x_1$,~\ldots,~$x_{n-i}\in H^2(BT_\alpha)\subset\m$
  clearly is $M_\m$-regular. Hence 
  \[
    \depth M_\m 
    \ge n-i.
  \]
  The description of depth via $\Ext$
  (\cf~the proof of Lemma~\ref{depth-cokernel})
  shows that the depth of a direct sum is the minimum of the depths
  of the summands. This gives the result.
\end{proof}

We now prove Claim~\ref{big-claim} by induction on~$i$, assuming that
$H^*_T(X,X_{i-1})_\m$ is zero or a {\CM} module of dimension~$d+n-i$.
For~$i=0$ this is true because we assume
$H^*_T(X,X_{-1})=H^*_T(X)$ to be free over the regular ring~$H^*(BT)$.

In the case~$H^*_T(X,X_{i-1})_\m=0$ the sequence~\eqref{exact-triple}
splits into short exact sequences for trivial reasons.
Suppose therefore that $H^*_T(X,X_{i-1})_\m$ is non-zero and Cohen-Macaulay
of dimension~$d+n-i$.
Since
\begin{equation}\label{bound-dim}
  \dim H^*_T(X,X_i)_\m\le\dim H^*_T(X,X_i)\le d+n-i-1
\end{equation}
by Lemma~\ref{first-claim} and $\dim H^*_T(X,X_{i-1})_\m= d+n-i$,
the image of~$H^*_T(X,X_i)_\m$ in~$H^*_T(X,X_{i-1})_\m$
has lower dimension than the target module,
hence is zero by Lemma~\ref{dimension-submodule}.
In other words, \eqref{short-sequence} is exact in this case as well.

Lemma~\ref{depth-cokernel} now gives
\begin{equation}\label{bound-depth}
  \depth H^*_T(X,X_i)_\m\ge d+n-i-1
\end{equation}
since $\depth H^*_T(X_i,X_{i-1})_\m\ge n-i$ by Lemma~\ref{second-claim}
and $\depth H^*_T(X,X_{i-1})_\m= d+n-i$ 
(or infinity if $H^*_T(X,X_{i-1})_\m=0$).

Comparing \eqref{bound-dim}~and~\eqref{bound-depth}
with the inequality~\eqref{depth-le-dim} finishes the proof.

\section{Actions without fixed points}

The Atiyah--Bredon sequence~\eqref{sequence-atiyah}
can never be exact
if $T$ acts without fixed points (unless $X=\emptyset$).
Yet there are interesting fixed-point-free actions.
In the context of compact differentiable $T$-manifolds and real coefficients,
Goertsches and Töben~\cite{GoertschesToeben:??}
observed that one can modify \eqref{sequence-atiyah}
to accommodate for this case. We brief\/ly comment
on such a generalization in our setting.

\begin{proposition}\label{dimension-orbit}
  Let $k$ be the minimal dimension of the $T$-orbits in~$X$.
  If $X$ satisfies condition~\eqref{condition-skeleton},
  then $\dim H_T^*(X) = d+n-k$.
\end{proposition}

\begin{proof}
  Since $X_{k-1}=\emptyset$, we have $\dim H_T^*(X) \le d+n-k$
  by Lemma~\ref{first-claim}.

  For the reverse inequality, 
  pick an~$x\in X_k$, say with isotropy group~$T'$.
  The map~$H^*(BT)\to H^*(BT')$ factors through the
  map~$H_T^*(X)\to H_T^*(T x)=H^*(BT')$ induced by the
  inclusion~$T x\hookrightarrow X$.
  This implies $\Ann H_T^*(X)\subset \Ann H^*(BT')$,
  hence $\dim H_T^*(X) \ge \dim H^*(BT')$.
  But, under condition~\eqref{condition-skeleton},
  Lemma~\ref{cohomological-dimension-isotropy} gives $\dim H^*(BT')=d+n-k$.
\end{proof}

A finitely generated graded module~$M$ over~$H^*(BT)$ is called \newterm{\CM}
if $M_\m$ is {\CM} over~$H^*(BT)_\m$ for all maximal ideals~$\m\in\Supp M$.
If $M$ is free over~$H^*(BT)$,
then it is {\CM} of maximal dimension~$d + n$
(\cf~the proof of Claim~\ref{big-claim}).
We remark in passing that the converse holds as well,
even if $R$ is not a field:
The Auslander--Buchsbaum formula~\cite[Thm~1.3.3]{BrunsHerzog:93}
together with~\cite[Cor.~IV.C.2]{Serre:65}
implies that $M$ 
is projective; freeness then follows from
the (comparatively easy) analogue
of the Quillen--Suslin theorem
for finitely generated \emph{graded} modules over polynomial rings
with coefficients in a PID \cite[Cor.~4.7]{Lam:78}.

\begin{theorem}
  \label{Goertsches-Toeben}
  Let $k$ be the minimal dimension of the $T$-orbits in~$X$.
  If $H_T^*(X)$ is {\CM} over~$H^*(BT)$
  and $X$ satisfies condition~\eqref{condition-skeleton},
  then the following sequence is exact:
  \begin{equation*}\label{sequence-goertsches-toeben}
  \let\longrightarrow\rightarrow
    0
    \longrightarrow H^*_T(X)
    \longrightarrow H^*_T(X_k)
    \longrightarrow H^{*+1}_T(X_{k+1}, X_k)
    \longrightarrow \cdots
    \longrightarrow H^{*+n}_T(X_n, X_{n-1})
    \longrightarrow 0.
  \end{equation*}
\end{theorem}

The only change to the proof of Theorem~\ref{Atiyah-Bredon}
is that the induction now starts at~$i=k$,
which is the dimension of~$H_T^*(X)$ by Proposition~\ref{dimension-orbit}.

\end{document}